\pdfoutput=1
\documentclass{IEEEtran}
\usepackage{setspace}
\ifCLASSINFOpdf
\else
\fi
\usepackage[cmex10]{amsmath}

 \usepackage[hyperindex=true,pdftitle={},
pdfauthor={St\'ephane Guerrier},colorlinks=TRUE,
pagebackref=false,citecolor=blue,plainpages=false,
pdfpagelabels]{hyperref} % colourlinks=false for printing

\usepackage{verbatim} % For code environment and math support
\usepackage{amsthm}
\usepackage{amsbsy}
\usepackage{amssymb}
\usepackage{amsbsy}
\usepackage{amsfonts}
\usepackage{graphics}
\usepackage[pdftex]{graphicx}
\usepackage{bm}
\usepackage[usenames]{color}
\usepackage{cite}
\usepackage{mathtools}
\usepackage{authblk}
\usepackage{booktabs}
\usepackage{nicefrac}
\usepackage{enumitem}
\usepackage{bbding}
\usepackage{soul}
\usepackage[numbers]{natbib}
\usepackage{nicefrac}
\usepackage{tikz}

\usepackage{xr}
\DeclareMathOperator*{\cov}{cov}

\DeclareMathOperator*{\var}{var}

\DeclareMathOperator{\AV}{AVar}
\DeclareSymbolFont{lettersA}{U}{txmia}{m}{it}
\DeclareMathSymbol{\real}{\mathord}{lettersA}{"92}
\DeclareMathSymbol{\field}{\mathord}{lettersA}{"83}

\def\boxit#1{\vbox{\hrule\hbox{\vrule\kern3pt
          \vbox{\kern3pt#1\kern3pt}\kern3pt\vrule}\hrule}}

\definecolor{pinegreen}{rgb}{0.0, 0.47, 0.44}	
	
\newtheoremstyle{mytheoremstyle} % name
    {0.3cm}                      % Space above
    {0cm}                        % Space below
    {\itshape}                   % Body font
    {}                           % Indent amount
    {\scshape}                   % Theorem head font
    {: }                          % Punctuation after theorem head
    {0em}                       % Space after theorem head
    {}  % Theorem head spec (can be left empty, meaning normal)

\theoremstyle{mytheoremstyle}

\newtheorem{Lemma}{Lemma}

\newtheoremstyle{myExampleRemarkstyle} % name
    {0.3cm}                    % Space above
    {0cm}                           % Space below
    {\itshape}                   % Body font
    {}                           % Indent amount
    {\scshape}                   % Theorem head font
    {: }                          % Punctuation after theorem head
    {0em}                       % Space after theorem head
    {}  % Theorem head spec (can be left empty, meaning normal)

\theoremstyle{myExampleRemarkstyle}

\newtheoremstyle{simuStyle}
{0.3cm} % Space above
{0cm} % Space below
{} % Body font
{} % Indent amount
{\bfseries} % Theorem head font
{.} % Punctuation after theorem head
{0em} % Space after theorem head
{} % Theorem head spec (can be left empty, meaning `normal')

\theoremstyle{simuStyle}

\makeatletter

\newenvironment{customLemma}[1]
  {\count@\c@Lemma
   \global\c@Lemma#1 %
    \global\advance\c@Lemma\m@ne
   \Lemma}
  {\endLemma
   \global\c@Lemma\count@}

\makeatother

\newtheoremstyle{stratStyle}
{0.3cm} % Space above
{0cm} % Space below
{} % Body font
{} % Indent amount
{\scshape} % Theorem head font
{: } % Punctuation after theorem head
{0em} % Space after theorem head
{} % Theorem head spec (can be left empty, meaning `normal')

\theoremstyle{stratStyle}

\DeclareSymbolFont{lettersA}{U}{txmia}{m}{it}
\DeclareMathSymbol{\real}{\mathord}{lettersA}{"92}
\DeclareMathSymbol{\field}{\mathord}{lettersA}{"83}

\begin{document}
\title{A Study of the Allan Variance for Constant-Mean Non-Stationary Processes}
%\title{Wavelet-based Approaches to overcome the Limits of Common Sensor Calibration Techniques}
%\title{Development of Efficient Sensor Calibration Techniques to overcome the Limits of Allan Variance-based Techniques}
%\title{Wavelet Variance Methods to respond to the Limits of Allan Variance-based Techniques for Sensor Calibration}
%\title{Limits of the Allan Variance Approach for Sensor Calibration and Development of Wavelet Variance based Estimators}
\author{Haotian Xu, St\'{e}phane Guerrier, Roberto Molinari \& Yuming Zhang  % <-this % stops a space
\thanks{\textbf{H. Xu} is a PhD student, Geneva School of Economics and Management, University of Geneva, 1211, Switzerland (E-mail: haotian.xu@unige.ch).}
\thanks{\textbf{S. Guerrier} is Assistant Professor, Department of Statistics, Pennsylvania State University, PA, 16801, USA.}
\thanks{\textbf{R. Molinari} is Visiting Assistant Professor, Department of Statistics \& Applied Probability, University of California, Santa Barbara, CA, 93117, USA.}
\thanks{\textbf{Y. Zhang} is a graduate student, Department of Statistics, Pennsylvania State University, PA, 16801, USA.}}
%\thanks{\textbf{M.-P. Victoria-Feser} is full Professor, Research Center for Statistics \& Geneva School of Economics and Management, University of Geneva, 1211, Switzerland (E-mail: Maria-Pia.VictoriaFeser@unige.ch).}}
%\thanks{\textbf{J. Skaloud} is Lecturer, Geodetic Engineering Laboratory, \'Ecole Polytechnique F\'ed\'erale de Lausanne, 1015, Switzerland (Email: jan.skaloud@epfl.ch).}}

%\thanks{Manuscript received .... ; revised ....}}


\markboth{}%
{Shell \MakeLowercase{\textit{et al.}}: Bare Demo of IEEEtran.cls for Journals}



\maketitle

% YS: je propose d'enlever le "The" et le "an".... mais c'est comme vous voulez, le titre tel quel est très bien...

\begin{abstract}
The Allan Variance (AV) is a widely used quantity in areas focusing on error measurement as well as in the general analysis of variance for autocorrelated processes in domains such as engineering and, more specifically, metrology. The form of this quantity is widely used to detect noise patterns and indications of stability within signals. However, the properties of this quantity are not known for commonly occurring processes whose covariance structure is non-stationary and, in these cases, an erroneous interpretation of the AV could lead to misleading conclusions. This paper generalizes the theoretical form of the AV to some non-stationary processes while at the same time being valid also for weakly stationary processes. Some simulation examples show how this new form can help to understand the processes for which the AV is able to distinguish these from the stationary cases and hence allow for a better interpretation of this quantity in applied cases.
\end{abstract}



% YS: added signal modeling

% craigmile2013estimation + moulines2008wavelet (long and short -memory AV/WV problems), abry1998long(GMWM), vernotte2012statistical (WN+RW)
\section{Introduction}
\label{Sec:intro}

The Allan Variance (AV) is a widely used quantity in areas going from engineering to physics where there is an interest in studying the stochastic stability of error measurements from various instruments such as, among others, clocks and oscillators. Its usefulness resides in the fact that it provides an extremely informative summary on the variance of time series or, more generally, of autocorrelated processes, especially when these are non-stationary and with infinite variance. Indeed, \cite{zhang2008allan} underlined how the AV is a better measure of uncertainty compared to standard methods (e.g. moving average variance) for processes such as random walks and non-stationary Fractional Autoregressive-Moving Average (ARFIMA) models, while being considerably useful also for stationary processes. For these processes the AV has a well known form which can help detect the kind of process, for example, from the log-log plot of the AV of an observed signal. The behaviour and forms of the AV for stationary and some non-stationary processes was studied in \cite{zhang2008allan} where the AV is used to detect and understand the process underlying a signal issued from different voltage measurements. However, there are many other applications for which the AV is of interest such as the detection of noise terms characterising inertial sensors \citep[see][]{el2008analysis} and many others \citep[see][for an overview]{percival2016wavelet}.

However, although the AV is extremely useful in the above settings, it is not known how it behaves when in the presence of other types of processes and whether it is able to distinguish between them. In this paper we intend to investigate the form of the AV for a particular class of processes which includes all those processes that have a constant mean but have a time-varying variance-covariance structure such as, for example, the Generalized AutoRegressive Conditional Heteroscedasticity (GARCH) models \citep[see][]{bollerslev1986generalized}, while processes with specific forms of mean or higher-order non-stationarity are not considered since they can either be dealt with through statistical regression techniques or simply cannot be detected by the AV. In particular, we focus on those processes which are characterized by a dependence structure by blocks since they are common in settings such as longitudinal studies or sensor calibration for navigation engineering. In the latter cases, the AV is often approximated by that of other known stationary processes such as, for example, the bias-instability process whose AV is often approximated by that of a first-order autoregressive process \citep[see, for example,][]{unsal2012estimation}. Moreover, it is not clear whether the AV can actually help to distinguish between these processes and those processes for which its form is currently known. The latter aspect is of particular relevance since it could lead to an erroneous interpretation of the observed process, for example assuming stationarity when this is not the case and reaching false conclusions.

In order to deal with the above mentioned non-stationary processes, this paper intends to study the theoretical form of the AV when the covariance structure is non-stationary. The consequent advantage of this study is that, by considering the varying covariance structure in the AV definition, it extends the applicability of those approaches that make use of the AV and raises awareness on its limitations, and inappropriate interpretation, in distinguishing and identifying these processes from stationary ones. With this in mind, Section \ref{sec.ovw.av} briefly defines the AV and describes its theoretical form for those processes which have been considered up to now. Section \ref{sec.ns.av} introduces the new theoretical form of both the overlapping and non-overlapping AV for processes whose covariance structure is non-stationary and shows how the form of the AV for stationary processes is a special case of this new form. In the same section, three case studies are presented which highlight the importance of these findings in order to better interpret processes through the AV. Finally Section \ref{sec.conc} concludes.

\section{Overview of the Allan Variance}
\label{sec.ovw.av}

To introduce the AV, let us first define $(X_t)_{t=1,\hdots,T}$ as a weakly stationary, discrete time and regularly spaced stochastic process with a constant mean $\mu$ (i.e. $\mathbb{E}[X_t] = \mu$) and an autocovariance function defined as follows
\begin{equation*}
    \gamma (h) = \cov (X_{t+h},X_t),
\end{equation*}
which depends solely on $h$, the distance between observations, with $\sigma_X^2 = \gamma (0)$ being the process variance. We can consequently define the autocorrelation function as
\begin{equation*}
    \rho(h) = \frac{\gamma (h)}{\gamma (0)}.
\end{equation*}
We consider the AV computed at dyadic scales ($\tau$) starting from local averages of the process which can be denoted as
\begin{equation}
    \bar{X}_{t}^{(n)} \equiv \frac{1}{n} \sum_{i = 1}^{n} X_{t - n + i}\, ,
    \label{mean.noav}
\end{equation}
where $n \in \left\{x \in \mathbb{N} \, : \;  1 \leq x < \log_2 (T) - 1 \right\}$ therefore determines the number of consecutive observations considered for the average. If the process has constant mean $\mu$, this implies that $\bar{X}_{t}^{(n)}$ also has the same mean and, based on these averages and following \cite{percival1994long}, the maximum-overlapping AV (MOAV) 
\begin{eqnarray}
        \AV_n \left(X_t \right) \equiv \frac{1}{2 m^{\star}} \sum_{k = 2n}^{T} \mathbb{E}\left[ \left(\bar{X}_{k}^{(n)} - \bar{X}_{k-n}^{(n)} \right)^2 \right].
    \label{eq:MOAVNS_def}
\end{eqnarray}
where $m^{\star} = T - 2n + 1$ and whose corresponding estimator is given by
\begin{eqnarray}
\label{eq:MOAVNS_est}
        \widehat{\AV}_n \left(X_t \right) = \frac{1}{2m^{\star}} \sum_{k = 2n}^{T} \left(\bar{x}_{k}^{(n)} - \bar{x}_{k-n}^{(n)} \right)^2.
\end{eqnarray}
where $\bar{x}_{t}^{(n)}$ denotes the sample equivalent of $\bar{X}_{t}^{(n)}$ based on a realization of the process $(X_t)$. 
Another version of the AV is the non-overlapping AV (NOAV) whose estimator however is not statistically as efficient as that of the MOAV (see App. \ref{ns.noav} for more details).

Keeping in mind the above definitions of the MOAV, \cite{zhang2008allan} delivered a general theoretical form of this quantity when applied to weakly stationary processes which is given by
\begin{eqnarray}
\label{stat.av}
   \AV_n \left(X_t \right)  &=& \frac{\sigma_X^2}{n^2} \bigg(n\left[1-\rho(n)\right] \\
   &+& \sum_{i=1}^{n-1} i \left[2 \rho(n-i) - \rho(i) - \rho(2n-i)\right]\bigg). \nonumber
\end{eqnarray}
Based on the above equation, the exact form of the AV for different stationary processes, such as the general class of AutoRegressive Moving Average (ARMA) models, can be derived. Moreover, \cite{zhang2008allan} provided the theoretical AV for non-stationary processes such as the random walk and ARFIMA models for which the AV, as mentioned earlier, represents a better measure of uncertainty compared to other methods.

Using the known theoretical forms of the AV, it is therefore possible to detect and distinguish different processes based on the pattern of their AV. Due to this, this quantity (or similar quantities such as the Haar wavelet variance) can be used as a mean to estimate different kind of processes \citep[see for example][]{guerrier2013wavelet}. However, there are many commonly encountered processes whose AV is not known and it is unclear to what point these can actually be distinguished from stationary processes. The next section delivers a more general form of the AV which includes these processes and studies if this quantity can actually be helpful in detecting them.

\section{Allan Variance for Constant-Mean Non-Stationary Processes}
\label{sec.ns.av}

As underlined in the previous sections, the AV is particularly useful for measuring uncertainty in non-stationary processes, especially when these have infinite variance. Nevertheless, there are other forms of non-stationarity for which the properties of the AV are unknown and these consist in those processes $(X_t)$ with a constant mean $\mu$ (independent of time) but a non-stationary covariance structure. This implies that the covariance function between observations at distance $h$ is also a function of time $t$ and can therefore be denoted as $\gamma(h,t)$. This type of process is very common in different areas going from engineering \citep[see][]{el2008analysis} to economics \citep[see][]{gallegati2014wavelet}.

To study the theoretical form of the AV for this class of processes, let us first define $\mathbf{X}_t^{(n)}$ as being the following vector of $n$ consecutive observations starting at $t-n +1$
\begin{equation*}
        \mathbf{X}_t^{(n)} \equiv \left[ X_{t-n+1} \; \cdots \; X_{t} \right]^T.
\end{equation*}
which contains the observations used to build the average in (\ref{mean.noav}). Using the above vector, for $t = n,..., T$, we can then define the matrix
\begin{eqnarray}
        \bm{\Sigma}_{t}^{(n)} \equiv \var \left( \mathbf{X}_{t}^{(n)} \right),
    \label{eq:mat:sigma}
\end{eqnarray}
and, for $t = 2n,..., T$, we define the matrix
\begin{eqnarray}
        \bm{\Gamma}_t^{(n)} \equiv \cov \left( \mathbf{X}_{t-n}^{(n)} , \,\mathbf{X}_{t}^{(n)} \right),
    \label{eq:mat:gamma}
\end{eqnarray}
These matrices represent the covariance matrices of the observations contained in each consecutive average. Indeed, $\bm{\Sigma}_{t}^{(n)}$ represents the covariance matrix of the observations within the average $\bar{X}_{t}^{(n)}$ which is used in definitions (\ref{eq:AVNS_def}) and (\ref{eq:MOAVNS_def}) while $\bm{\Gamma}_t^{(n)}$ represents the cross-covariance between these two sets of observations. A visual representation of these quantities is given in App. \ref{app.graph}.

In this section we will only consider the non-stationary MOAV while the form of the non-stationary NOAV is given in App. \ref{ns.noav}. Based on the above matrices, we can also define different quantities according to the matrix of reference and the lags between observations. More specifically, let us first consider the case in which we are interested in lags $h$ such that $0 \leq h < n$. Because of the overlapping nature of the AV, the observations at these lags can belong to the sets of observations within both the matrix $\bm{\Sigma}_{t}^{(n)}$ and  the matrix $\bm{\Gamma}_t^{(n)}$ and, for these sets of observations within the matrix $\bm{\Sigma}_{t}^{(n)}$, we can define the following quantity
\begin{eqnarray*}
    \widetilde{\gamma}(h) \equiv \frac{1}{2m^{\star}(n-h)} \sum_{t = 2n}^{T} \sum_{s=0}^{n-h-1}\cov (X_{t-n-s-h}, \,X_{t-n-s})\\
     + \cov (X_{t-s-h}, \,X_{t-s}).
    \label{eq:gamma_avg_4}
\end{eqnarray*}
If, however, the observations at the considered lags are among the set of observations only within the matrix $\bm{\Gamma}_t^{(n)}$, we define the quantity below
\begin{equation*}
    \widetilde{\gamma}^{\ast}(h) \equiv \frac{1}{m^{\star}h} \sum_{t = 2n}^{T} \sum_{s=1}^{h}\cov \left(X_{t-n+s-h}, \,X_{t-n+s} \right).
    \label{eq:gamma_avg_5}
\end{equation*}
Finally, when considering lags $h$ such that $n \leq h \leq 2n-1$, the set of observations at these lags can only be considered within the matrix $\bm{\Gamma}_t^{(n)}$ and, for this final case, we define the quantity
\begin{equation*}
    \widetilde{\gamma}(h) \equiv \frac{1}{m^{\star}(2n-h)} \sum_{t = 2n}^{T} \sum_{s=0}^{2n-h-1}\cov \left(X_{t-s-h}, \,X_{t-s}\right).
    \label{eq:gamma_avg_6}
\end{equation*}
The above definitions can be seen as generalized definitions of the autocovariance which consists in the average autocovariance for a given lag $h$. Because of this, it must be underlined that these definitions are not at all equivalent to the covariance function $\gamma(h,t)$ but correspond to an average of this function over all times $t$. Having specified this, we can now provide the following lemma.
\begin{Lemma}
\label{lem:AVNS_MO}
  The non-stationary MOAV is given by
  \begin{eqnarray*}
      \AV_n =& \frac{1}{2 m^{\star} n^2} \Bigg\{ 2nm^{\star} \, \widetilde{\gamma}(0) + 2\Bigg[\sum_{h=1}^{n-1}2m^{\star}(n-h) \, \widetilde{\gamma}(h) \\
      &- m^{\star}h \, \widetilde{\gamma}^{\ast}(h)
      - \sum_{h=n}^{2n-1} m^{\star}(2n-h) \, \widetilde{\gamma}(h) \Bigg] \Bigg\}.
  \end{eqnarray*}
\end{Lemma}
The proof of this lemma is given in App. \ref{app.moavns}. Considering this expression, an aspect that must be underlined is that the definitions of the functions $\widetilde{\gamma}(h)$ and $\widetilde{\gamma}^{\ast}(h)$ given earlier simplify to the autocovariance function $\gamma(h)$ when dealing with a weakly stationary process. In the latter case, the form of the non-stationary AV consequently reduces to the expression in (\ref{stat.av}) for which a more detailed discussion is given in App. \ref{app.statequiv}. As a final note to this result it should also be highlighted that, in some of the considered non-stationary cases, the estimators of the MOAV defined in (\ref{eq:MOAVNS_est}) and of the NOAV (see App. \ref{ns.noav}) do not necessarily have the same expectation. 

Having underlined these points and with the general definition of the AV given in Lemma \ref{lem:AVNS_MO}, we can now investigate its properties when assuming the process of interest is within the class of processes treated in this paper. The next sections report some simulation studies regarding some of these processes in attempt to understand also whether the AV is a useful quantity to detect them and distinguish them from weakly stationary processes. In all cases, the simulated process is of length $T = 1,000$ (except for the bias-instability process where $T = 2,500$) and is simulated 50 times. The estimated AV is represented in plots along with the theoretical stationary and non-stationary forms of the AV in order to understand its behaviour under these different assumptions. 

\subsubsection{Non-Stationary White Noise}

\begin{figure}
  \centering
    \includegraphics[width=0.45\textwidth]{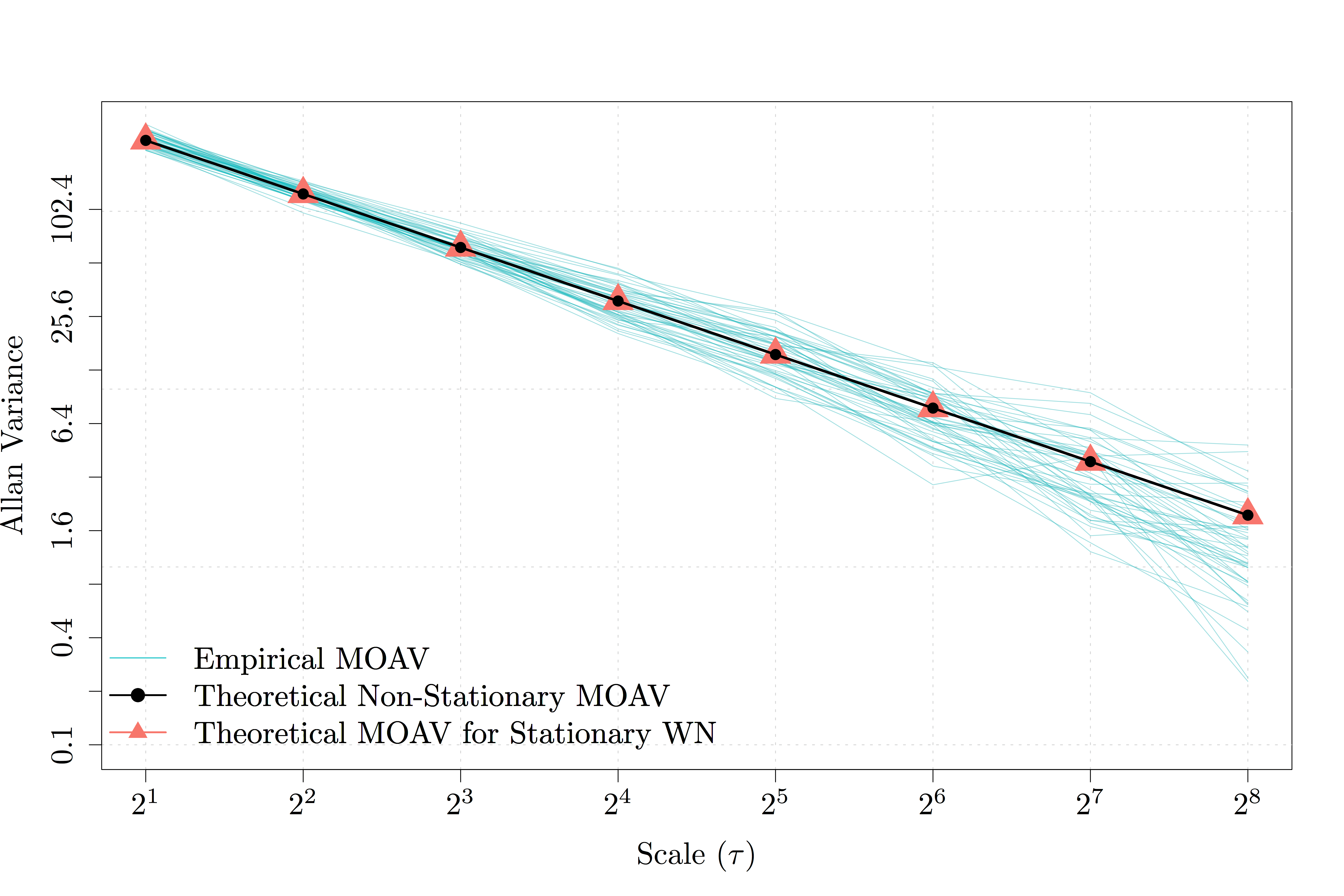}
      \caption{Logarithm of the MOAV of the non-stationary white noise process for scales $\tau = 2^n$. Estimated MOAV (light-blue lines); theoretical non-stationary MOAV (black line with dots) and theoretical stationary AV based on the average variance (red line with triangles).}
      \label{fig.wn}
\end{figure}

The first process we study is the non-stationary white noise by which we intend, without loss of generality, all those zero-mean processes whose variance changes with time. The evolution of the variance in time can either be completely random or can follow a specific parametric model such as, for example, a GARCH process. The goal of studying the AV for these types of processes is to understand whether it is able to detect such a structure in a time series and if it can be distinguished from a stationary white noise process. For this purpose, the true non-stationary process considered in the simulation study is generated from the following model
$$X_t \sim \mathcal{N}(0,\sigma_t^2),$$
where $\sigma_t^2 = t$, with $t = 1, \hdots, T$. The theoretical stationary form is based on the average of the variances used to simulate the processes (i.e. $\bar{\sigma}^2 = \nicefrac{(1+T)}{2}$ in this example). Fig. \ref{fig.wn} represents the estimated AVs along with the theoretical forms (stationary and non-stationary).
In this case, it can be seen how both theoretical forms correspond and the estimated AVs closely follow these quantities. This example confirms that the AV is therefore unable to distinguish between a stationary white noise process and a white noise process whose second-order behaviour is non-stationary.

\subsubsection{Bias-Instability}

\begin{figure}
  \centering
    \includegraphics[width=0.45\textwidth]{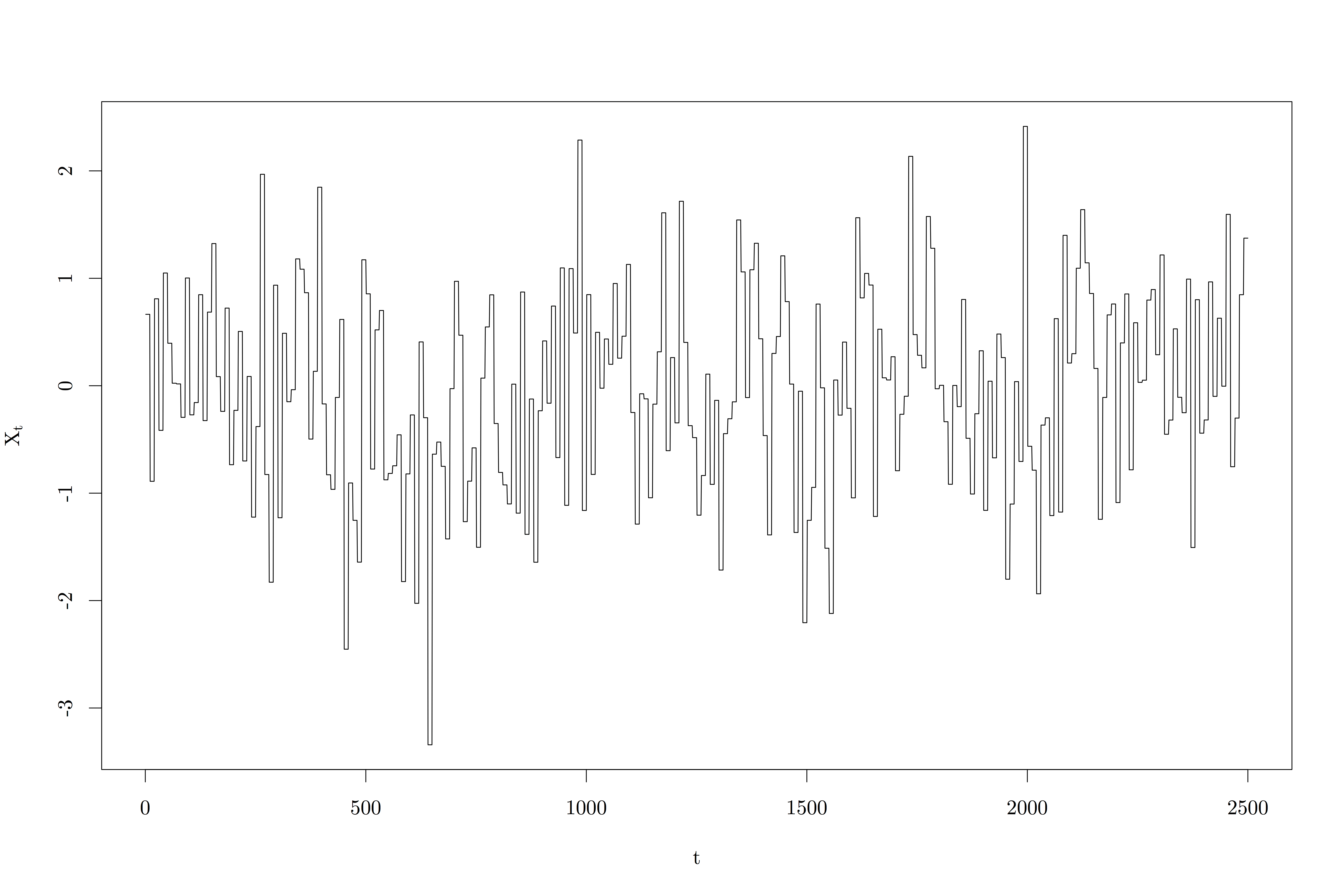}\vspace{-0.2cm}
    \includegraphics[width=0.45\textwidth]{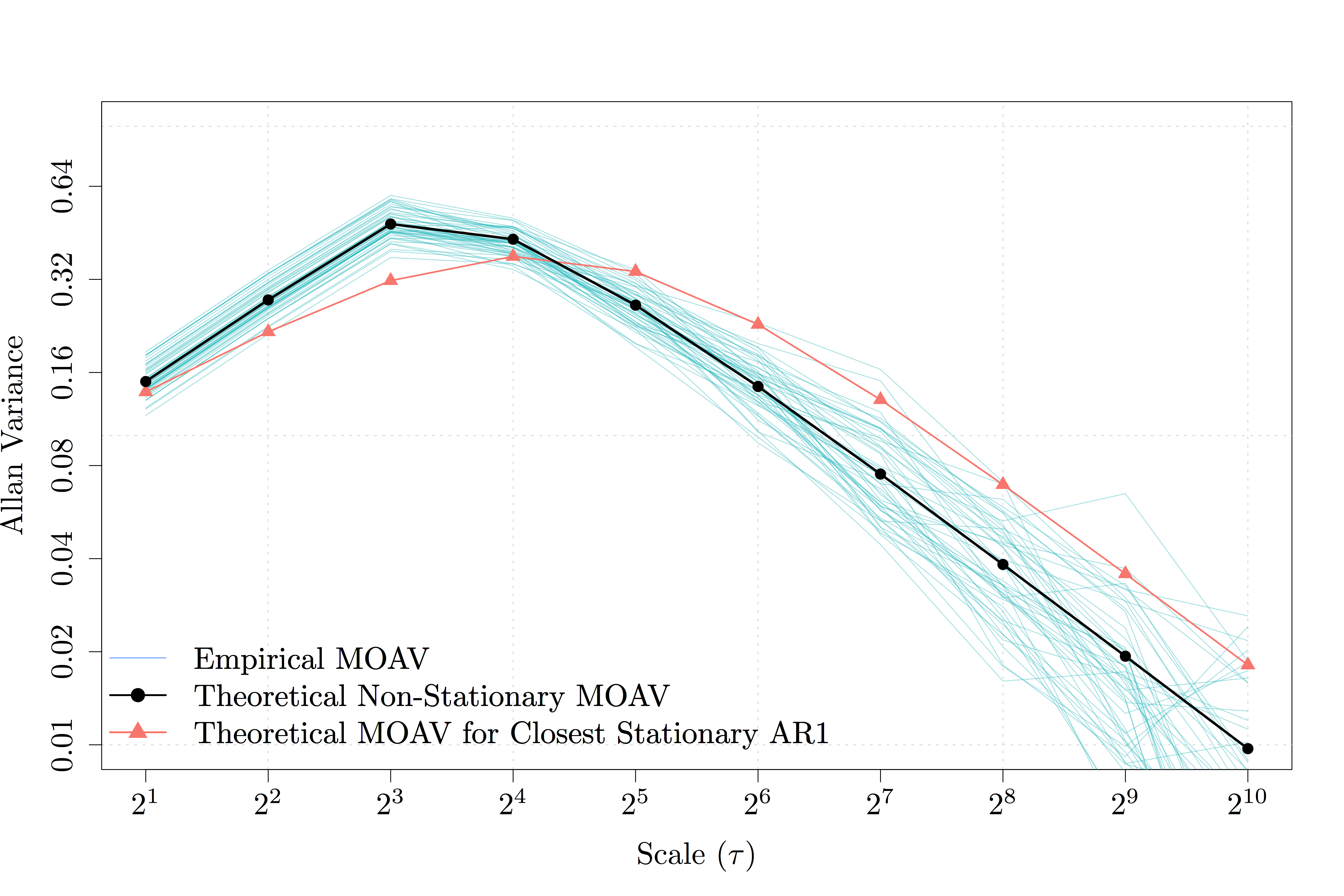}
      \caption{Top: Realization of the bias-instability process with $\sigma^2 = 1$ and the length of block $b_i = 10$, $\forall i = 1, \hdots, B$, and $B = 250$. Bottom: Logarithm of the MOAV of the bias-instability process for scales $\tau = 2^n$. Estimated MOAV (light-blue lines); theoretical non-stationary MOAV (black line with dots) and theoretical stationary MOAV of an AR1 approximating the bias-instability MOAV (red line with triangles).}
      \label{fig:bi}
\end{figure}

The bias-instability process is a commonly known process in the engineering domain, specifically for inertial sensor calibration for navigation. The characteristic of this process is that it consists in different concatenated sequences (blocks) where, within each block, the realization of a random variable is repeated (i.e. constant). More formally, let $b_i, \,\, i = 1, \hdots, B$, represent the set of time indices belong to the $i^{th}$ block within the time series, and let $C_i \overset{iid}{\sim} \mathcal{N}(0,\sigma^2)$. We can then define this process as
\[   
X_t = c_i \quad \text{if} \,\, t \in b_i.
\]
One realization of the bias-instability process is illustrated in the top panel Fig. \ref{fig:bi} where the length of block $b_i$ is 10 for all $i = 1, \hdots, B$, and $B = 250$.
Since the theoretical form of the AV for this process is not known exactly, it is often approximated by the AV of a First-Order AutoRegressive (AR1) process. Although this approximation can be useful, it is nevertheless still an approximation and, using the form given in Lemma \ref{lem:AVNS_MO}, we can now obtain a theoretical form for the AV of this process which is represented in the bottom panel of Fig. \ref{fig:bi}. Indeed the latter plot shows that the estimated AVs closely follow the theoretical non-stationary form given earlier. The red line represents the AV of a stationary AR1 process which is supposed to approximate the true AV of bias-instability. The latter is the result of the averaging of the theoretical AV for a stationary AR1 process estimated via maximum-likelihood on each of the simulated processes. It is clear how, although close over some scales, this approximation is not good enough when considering the logarithmic representation of the AV. Therefore, knowing the exact form of the AV for this process would allow to better interpret the signals characterised by bias-instability. 

\subsubsection{Block-Structure Autoregressive Processes}

\begin{figure}
  \centering
    \includegraphics[width=0.45\textwidth]{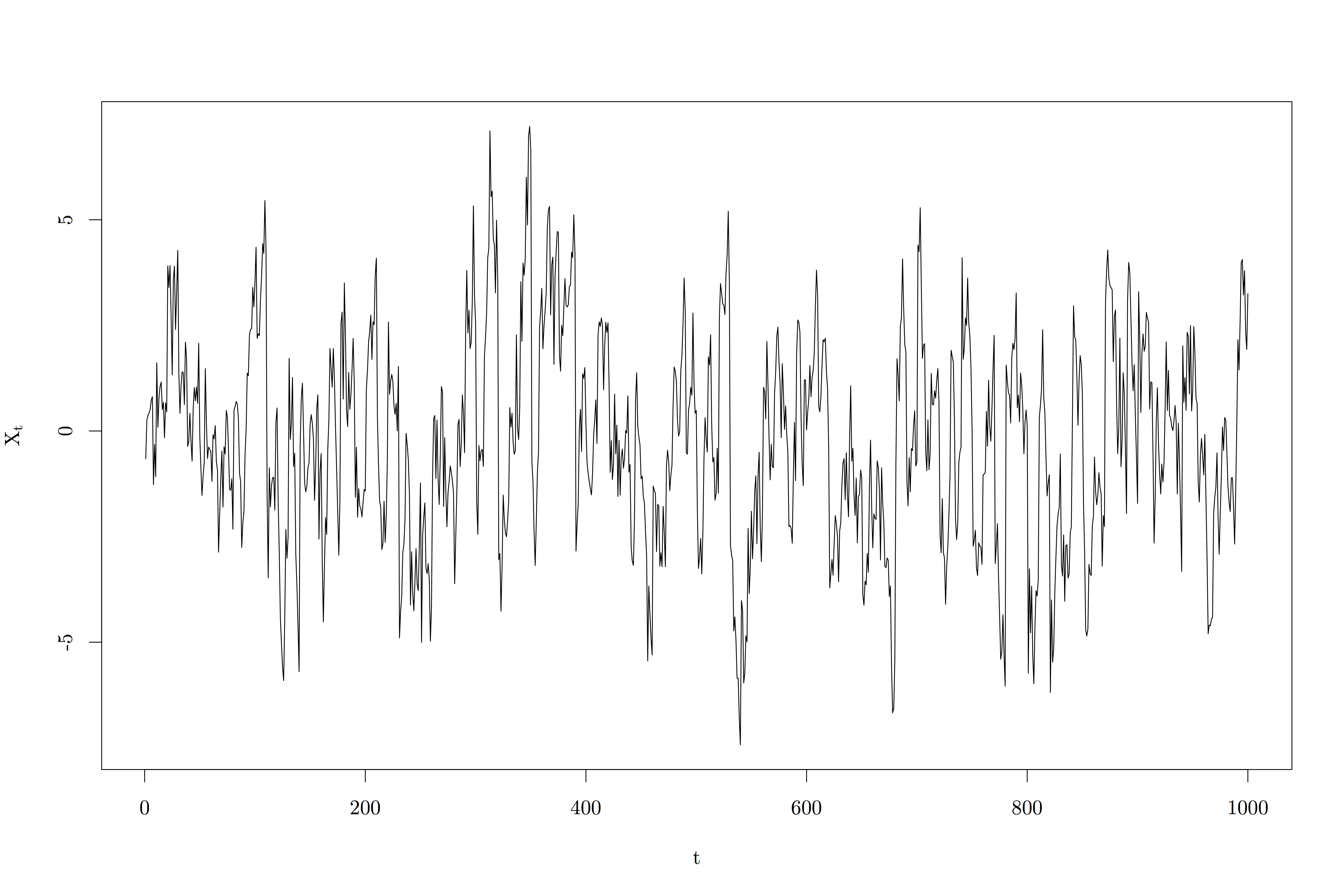}\vspace{-0.2cm}
    \includegraphics[width=0.45\textwidth]{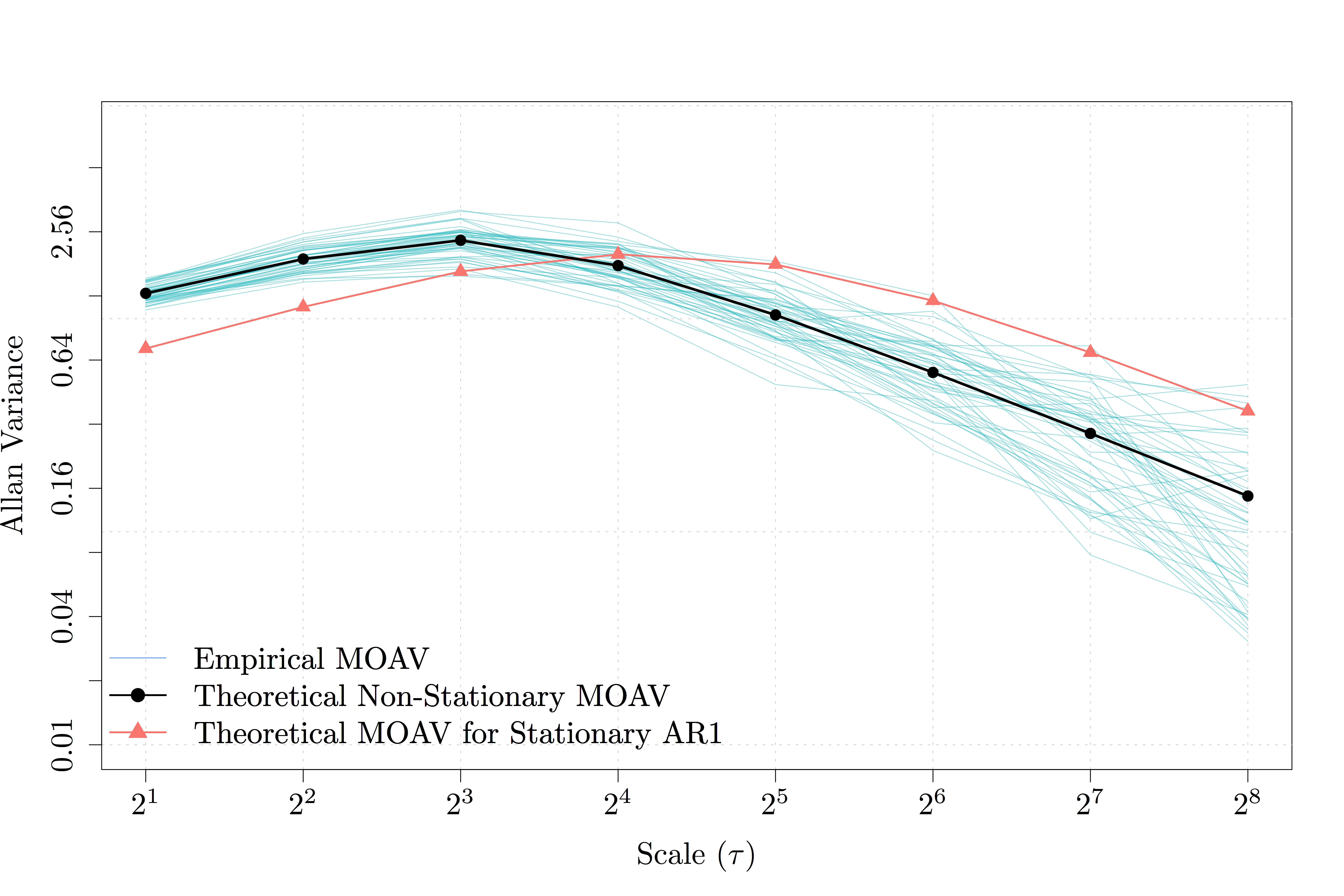}
      \caption{Top: Realization of the block-structure autoregressive process with $\phi = 0.9$, $\sigma^2 = 1$ and the length of block $b_i = 10, \forall i = 1, \hdots, B$, and $B = 100$. Bottom: Logarithm of the MOAV of the block-structure first-order autoregressive process  for scales $\tau = 2^n$. Estimated MOAV (light-blue lines); theoretical non-stationary MOAV (black line with dots) and theoretical stationary MOAV assuming no block structure (red line with triangles).}
      \label{fig:ar1}
\end{figure}

As a final example we consider a block-structure AR1 process. Similarly to the bias-instability process, within this paper, we define a block-structure process as a process whose parameters are fixed but is made by concatenated time periods (blocks) where observations within each block are generated independently from those in the other blocks. An example is given by the settings of longitudinal studies in which each subject can be measured over time and, although the subjects are independent from each other, these measurements can be explained by an autocorrelated process within each subject. To define this process formally, let $X_t^{(i)} \sim F_{\bm{\theta}}$ denote the following AR1 process
$$X_t^{(i)} = \phi X_{t-1}^{(i)} + \epsilon_t$$
with parameter vector $\bm{\theta} = [\phi \,\,\,\, \sigma^2 ]^T$ where $\phi \in (-1,1)$ and $\epsilon_t \overset{iid}{\sim} \mathcal{N}(0,\sigma^2)$. If again we let $b_i$ denote the $i^{th}$ block, then the block-structure AR1 process can be defined as
\[   
X_t = X_t^{(i)} \quad \text{if} \,\, t \in b_i.
\]{}
where $X_t^{(i)}$ is independent of $X_t^{(j)}$, $\forall \, i \neq j$. By defining $\phi = 0.9$ and $\sigma^2 = 1$ for the simulation study, the top panel of Fig. \ref{fig:ar1} shows a realization of this process while the bottom panel of Fig. \ref{fig:ar1} illustrates the results of the simulations for this particular process.
As for bias-instability, it can be observed how the stationary form of the AV (that does not consider the block structure) is not close to the estimated AVs while the non-stationary form provided in this paper adequately represents this process and can therefore allow to distinguish between a stationary autoregressive process and a block-structure one.

\section{Conclusions}
\label{sec.conc}

Within this paper we wanted to underline an issue concerning the AV which had not yet been studied. Indeed, the behaviour of the AV in commonly occurring settings where the covariance structure of the processes is non-stationary was unknown and, in many cases, was either ignored or dealt with through approximations. The consequence of the latter approaches would probably consist in erroneous interpretations and conclusions drawn from an AV analysis. For this reason, this paper studied the form of the AV for this class of processes thereby generalizing its form also for weakly stationary processes. Based on this, several examples were provided in which the properties of the AV were studied, highlighting its ability to detect these processes and to eventually distinguish them from stationary ones, making researchers and practitioners more aware of issues related to the interpretation and use of this quantity in more general and common settings.

\bibliographystyle{spmpsci}
\bibliography{ref}

\appendices
%%%%%%%%%%%%%%%%%%%%%%%%%%%%%%%%%%%%%%%%%%%%%%%%%%%%%%%%%%%%%
%%%%%%%%%%%%%%%%%%%%%%%%%%%%%%%%%%%%%%%%%%%%%%%%%%%%%%%%%%%%%

%%% THIS COMMANDE MAKES NUMERING EQUATION (B-x)
  \renewcommand{\theequation}{D-\arabic{equation}}
  % redefine the command that creates the equation no.
  \setcounter{equation}{0}  % reset counter

%%%%%%%%%%%%%%%%%%%%%%%%%%%%%%%%%%%%%%%%%%%%%%%%%%%%%%%%%%%%%
%%%%%%%%%%%%%%%%%%%% PROOF COROLLARY 1 %%%%%%%%%%%%%%%%%%%%%%
%%%%%%%%%%%%%%%%%%%%%%%%%%%%%%%%%%%%%%%%%%%%%%%%%%%%%%%%%%%%%

\section{Graphical illustration of MOAV}
\label{app.graph}

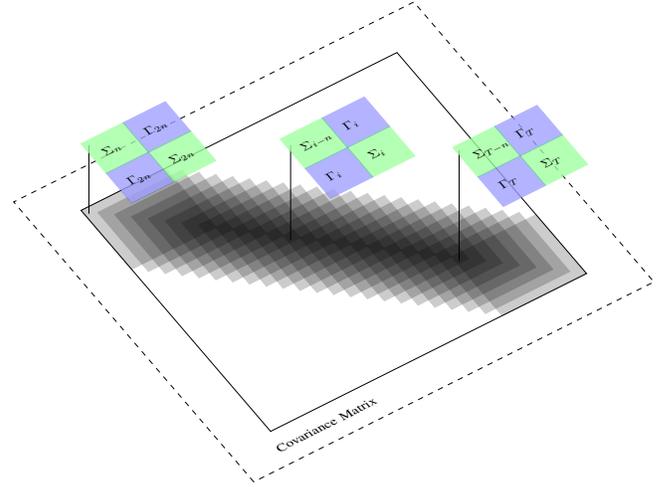
\begin{figure}[!h]
     \scalebox{0.7}{\newcommand{\Covariance}{Covariance Matrix}

\newcommand{\yslant}{0.5}
\newcommand{\xslant}{-0.6}

\centering
\begin{tikzpicture}[scale=0.2]
    \pgfmathsetmacro{\n}{8}

	% Cross Product Level
	\begin{scope}[
		yshift=-270,
		every node/.append style={yslant=\yslant,xslant=\xslant},
		yslant=\yslant,xslant=\xslant
	] 
		% The frame:
		\draw[black, dashed, thin] (-6,-5) rectangle (32,33); 
		\draw[black, thin] (-2,-1) rectangle (28,29);

        % function
        \begin{scope}
            \foreach \x in {0, ..., 22}
                \filldraw[thin,black,opacity=.2] (\x-2, 21-\x) rectangle (\x-2+\n, 21-\x+\n);

        \end{scope}
        
        % Agents:
		%\draw[fill=red]  
		%	(-1.6,28.5) circle (.2) ;
        
        % Labels:
	    \fill[black]
		    (-3,-3) node[right, scale=.7] {\Covariance};
	\end{scope}
	
    % vertical line for linking agents on the 2 levels
	\draw[ultra thin](-18.65,9.5) to (-18.65,16);
	
	% Square level
	\begin{scope}[
		yshift=-90,
		every node/.append style={yslant=\yslant,xslant=\xslant},
		yslant=\yslant,xslant=\xslant
	]
		% The frame:
		\fill[white,fill opacity=.7] (-2,21) rectangle (6,29);  % Opacity
		\filldraw[thin,blue,opacity=.3] (-2,21) rectangle (2,25);
		\filldraw[thin,blue,opacity=.3] (2,25) rectangle (6,29);
        \filldraw[thin,green,opacity=.3] (-2,25) rectangle (2,29);
		\filldraw[thin,green,opacity=.3] (2,21) rectangle (6,25);

        % Agents:
		%\draw[fill=red]  
		%	(-1.6,28.5) circle (.2) ;
	\fill[black]
		    (2.5,27) node[right, scale=.7] {$\Gamma_{2n}$};
    \fill[black]
            (-1.5,23) node[right, scale=.7] {$\Gamma_{2n}$};
    \fill[black]
            (-1.5,27) node[right, scale=.7] {$\Sigma_{n}$};
    \fill[black]
            (2.5,23) node[right, scale=.7] {$\Sigma_{2n}$};
		
	\end{scope} 
	
	    % vertical line for linking agents on the 2 levels
	\draw[ultra thin](0.5,7) to (0.5,16);
	
	% Square level
	\begin{scope}[
		yshift=-90,
		every node/.append style={yslant=\yslant,xslant=\xslant},
		yslant=\yslant,xslant=\xslant
	]
		% The frame:
		\fill[white,fill opacity=.7] (11.5,12) rectangle (19.5,20);  % Opacity
		\filldraw[thin,blue,opacity=.3] (11.5,12) rectangle (15.5,16);
		\filldraw[thin,blue,opacity=.3] (15.5,16) rectangle (19.5,20);
        \filldraw[thin,green,opacity=.3] (11.5,16) rectangle (15.5,20);
		\filldraw[thin,green,opacity=.3] (15.5,12) rectangle (19.5,16);

        % Agents:
		%\draw[fill=red]  
		%	(-1.6,28.5) circle (.2) ;
	\fill[black]
		    (16,14) node[right, scale=.7] {$\Sigma_{i}$};
    \fill[black]
            (12,18) node[right, scale=.7] {$\Sigma_{i-n}$};
    \fill[black]
            (16,18) node[right, scale=.7] {$\Gamma_i$};
    \fill[black]
            (12,14) node[right, scale=.7] {$\Gamma_i$};
		
	\end{scope} 
	
	% vertical line for linking agents on the 2 levels
	\draw[ultra thin](16.5,5) to (16.5,15.5);
	
	% Square level
	\begin{scope}[
		yshift=-90,
		every node/.append style={yslant=\yslant,xslant=\xslant},
		yslant=\yslant,xslant=\xslant
	]
		% The frame:
		\fill[white,fill opacity=.7] (22.5,3) rectangle (30.5,11);  % Opacity
		\filldraw[thin,blue,opacity=.3] (22.5,3) rectangle (26.5,7);
		\filldraw[thin,blue,opacity=.3] (26.5,7) rectangle (30.5,11);
        \filldraw[thin,green,opacity=.3] (22.5,7) rectangle (26.5,11);
		\filldraw[thin,green,opacity=.3] (26.5,3) rectangle (30.5,7);

        % Agents:
		%\draw[fill=red]  
		%	(-6.6,33.5) circle (.2) ;
	\fill[black]
		    (27,5) node[right, scale=.7] {$\Sigma_{T}$};
    \fill[black]
            (23,9) node[right, scale=.7] {$\Sigma_{T-n}$};
    \fill[black]
            (23,5) node[right, scale=.7] {$\Gamma_T$};
    \fill[black]
            (27,9) node[right, scale=.7] {$\Gamma_T$};
		
	\end{scope} 
	
\end{tikzpicture}}  

    \caption{Graphical illustration of matrices $\bm{\Sigma}_{t}^{(n)}$ and $\bm{\Gamma}_t^{(n)}$ for the MOAV.}
    \label{fig1}
\end{figure}

To graphically illustrate the quantities defined in Section \ref{sec.ns.av}, Fig. \ref{fig1} represents the true covariance matrix for a given process and highlights how the AV is related to this matrix by overlapping square matrices along the diagonal, each of which is composed of the quantities defined in Eq. (\ref{eq:mat:sigma}) and Eq. (\ref{eq:mat:gamma}).

\section{Theoretical form of the NOAV for Non-Stationary Processes}
\label{ns.noav}

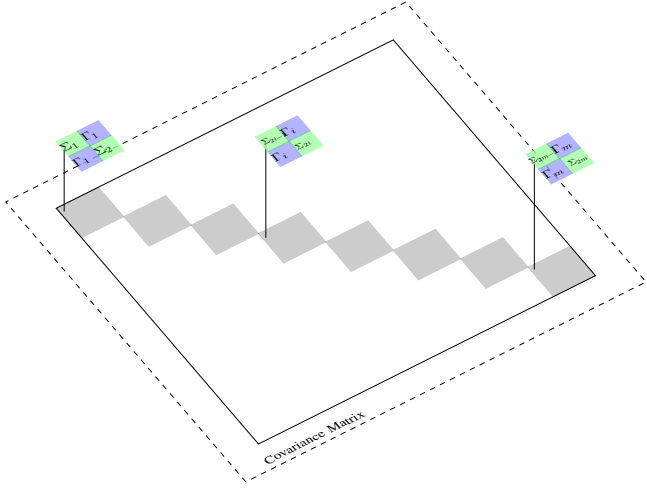
\begin{figure}[!h]
     \scalebox{0.7}{\newcommand{\Covariance}{Covariance Matrix}

\newcommand{\yslant}{0.5}
\newcommand{\xslant}{-0.6}

\centering
\begin{tikzpicture}[scale=0.2]
    \pgfmathsetmacro{\n}{4}

	% Cross Product Level
	\begin{scope}[
		yshift=-270,
		every node/.append style={yslant=\yslant,xslant=\xslant},
		yslant=\yslant,xslant=\xslant
	] 
		% The frame:
		\draw[black, dashed, thin] (-5,-6) rectangle (33,32); 
		\draw[black, thin] (-2,-3) rectangle (30,29);

        % function
        \begin{scope}
            \foreach \x in {0, 4,..., 30}
                \filldraw[thin,black,opacity=.2] (\x-2, 25-\x) rectangle (\x-2+\n, 25-\x+\n);

        \end{scope}
        
        % Agents:
		%\draw[fill=red]  
		%	(-1.6,28.5) circle (.2) ;
        
        % Labels:
	    \fill[black]
		    (-3,-5) node[right, scale=.7] {\Covariance};
	\end{scope}
	
    % vertical line for linking agents on the 2 levels
	\draw[ultra thin](-18.65,9.5) to (-18.65,16);
	
	% Square level
	\begin{scope}[
		yshift=-90,
		every node/.append style={yslant=\yslant,xslant=\xslant},
		yslant=\yslant,xslant=\xslant
	]
		% The frame:
		\fill[white,fill opacity=.7] (-2,25) rectangle (2,29);  % Opacity
		\filldraw[thin,blue,opacity=.3] (0,27) rectangle (2,29);
		\filldraw[thin,blue,opacity=.3] (-2,25) rectangle (0,27);
        \filldraw[thin,green,opacity=.3] (-2,27) rectangle (0,29);
		\filldraw[thin,green,opacity=.3] (0,25) rectangle (2,27);

        % Agents:
		%\draw[fill=red]  
		%	(-1.6,28.5) circle (.2) ;
	\fill[black]
		    (-0.5,28) node[right, scale=.7] {$\Gamma_{1}$};
    \fill[black]
            (-2.5,26) node[right, scale=.7] {$\Gamma_{1}$};
    \fill[black]
            (-2.5,28) node[right, scale=.7] {$\Sigma_{1}$};
    \fill[black]
            (-0.5,26) node[right, scale=.7] {$\Sigma_{2}$};
		
	\end{scope} 
	
	    % vertical line for linking agents on the 2 levels
	\draw[ultra thin](0.5,7) to (0.5,16);
	
	% Square level
	\begin{scope}[
		yshift=-90,
		every node/.append style={yslant=\yslant,xslant=\xslant},
		yslant=\yslant,xslant=\xslant
	]
		% The frame:
		\fill[white,fill opacity=.7] (11.5,16) rectangle (15.5,20);  % Opacity
		\filldraw[thin,blue,opacity=.3] (13.5,18) rectangle (15.5,20);
		\filldraw[thin,blue,opacity=.3] (11.5,16) rectangle (13.5,18);
        \filldraw[thin,green,opacity=.3] (11.5,18) rectangle (13.5,20);
		\filldraw[thin,green,opacity=.3] (13.5,16) rectangle (15.5,18);

        % Agents:
		%\draw[fill=red]  
		%	(-1.6,28.5) circle (.2) ;
	\fill[black]
		    (13.3,17) node[right, scale=.5] {$\Sigma_{2i}$};
    \fill[black]
            (11.3,19) node[right, scale=.5] {$\Sigma_{2i-1}$};
    \fill[black]
            (13,19) node[right, scale=.7] {$\Gamma_i$};
    \fill[black]
            (11,17) node[right, scale=.7] {$\Gamma_i$};
		
	\end{scope} 
	
	% vertical line for linking agents on the 2 levels
	\draw[ultra thin](26,4) to (26,14.5);
	
	% Square level
	\begin{scope}[
		yshift=-90,
		every node/.append style={yslant=\yslant,xslant=\xslant},
		yslant=\yslant,xslant=\xslant
	]
		% The frame:
		\fill[white,fill opacity=.7] (28.5,1.5) rectangle (32.5,5.5);  % Opacity
		\filldraw[thin,blue,opacity=.3] (30.5,3.5) rectangle (32.5,5.5);
		\filldraw[thin,blue,opacity=.3] (28.5,1.5) rectangle (30.5,3.5);
        \filldraw[thin,green,opacity=.3] (28.5,3.5) rectangle (30.5,5.5);
		\filldraw[thin,green,opacity=.3] (30.5,1.5) rectangle (32.5,3.5);

        % Agents:
		%\draw[fill=red]  
		%	(-1.6,28.5) circle (.2) ;
	\fill[black]
		    (30.5,2.5) node[right, scale=.5] {$\Sigma_{2m}$};
    \fill[black]
            (28,4.5) node[right, scale=.5] {$\Sigma_{2m-1}$};
    \fill[black]
            (28,2.5) node[right, scale=.7] {$\Gamma_m$};
    \fill[black]
            (30,4.5) node[right, scale=.7] {$\Gamma_m$};
		
	\end{scope} 
	
\end{tikzpicture}}  

    \caption{Graphical illustration of matrices $\bm{\Sigma}_{2k}^{(n)}$, $\bm{\Sigma}_{2k-1}^{(n)}$ and $\bm{\Gamma}_k^{(n)}$ for the NOAV.}
    \label{fig2}
\end{figure}

The non-overlapping AV (NOAV) is defined as:
\begin{eqnarray}
    \AV_n \left(X_t \right) \equiv \frac{1}{2m} \sum_{k = 1}^{m} \mathbb{E} \left[ \left(\bar{X}_{2k}^{(n)} - \bar{X}_{2k-1}^{(n)} \right)^2\right],
    \label{eq:AVNS_def}
\end{eqnarray}
where $m = \lfloor\frac{T}{2n}\rfloor$. The corresponding estimator for this quantity is given by
\begin{eqnarray*}
    \widetilde{\AV}_n \left(X_t \right) = \frac{1}{2m} \sum_{k = 1}^{m} \left(\bar{x}_{2k}^{(n)} - \bar{x}_{2k-1}^{(n)} \right)^2.
    \label{eq:AVNS_est}
\end{eqnarray*}
This estimator is less efficient than the MOAV, mainly because it is based on fewer averages and therefore on a smaller sample size. To define the theoretical form of the NOAV for the non-stationary processes of interest, we first define the vector  $\mathbf{X}_j^{(n)}$ of $n$ consecutive observations starting at $(j-1)n +1$, i.e.
\begin{equation*}
        \mathbf{X}_j^{(n)} \equiv \left[ X_{(j-1) n + 1} \; \cdots \; X_{j n} \right]^T.
\end{equation*}
Using the above, for $k = 1,..., m$, we define the matrices $\bm{\Sigma}_{2k}^{(n)}$, $\bm{\Sigma}_{2k-1}^{(n)}$ and $\bm{\Gamma}_k^{(n)}$ as follows:
\begin{eqnarray*}
        \bm{\Sigma}_{2k}^{(n)} \equiv \var \left( \mathbf{X}_{2k}^{(n)} \right), \;\; \bm{\Sigma}_{2k-1}^{(n)} \equiv \var \left( \mathbf{X}_{2k-1}^{(n)} \right)\;\; \\ \text{and} \;\;
        \bm{\Gamma}_k^{(n)} \equiv \cov \left( \mathbf{X}_{2k-1}^{(n)} , \mathbf{X}_{2k}^{(n)} \right).
    \label{eq:mat:sigma:gamma}
\end{eqnarray*}
As in Appendix \ref{app.graph}, the above matrices are graphically represented in Fig. \ref{fig2} where, as opposed to the MOAV, these matrices do not overlap along the diagonal of the covariance matrix of the process.
We then let $\bar{\sigma}^{(n)}_{2k}$, $\bar{\sigma}^{(n)}_{2k-1}$ and $\bar{\gamma}^{(n)}_k$ denote the averages of the matrices $\bm{\Sigma}_{2k}^{(n)}$, $\bm{\Sigma}_{2k-1}^{(n)}$ and $\bm{\Gamma}_k^{(n)}$, respectively, i.e.
\begin{eqnarray*}
    \bar{\sigma}^{(n)}_{2k} \equiv \frac{1}{n^2}\sum_{i = 1}^{n} \sum_{j = 1}^{n} \, \left( \bm{\Sigma}_{2k}^{(n)} \right)_{i,j},  \;\; \\ \bar{\sigma}^{(n)}_{2k-1} \equiv \frac{1}{n^2}\sum_{i = 1}^{n} \sum_{j = 1}^{n} \, \left( \bm{\Sigma}_{2k-1}^{(n)} \right)_{i,j}, \;\;\\
    \bar{\gamma}^{(n)}_k \equiv \frac{1}{n^2} \sum_{i = 1}^{n} \sum_{j = 1}^{n} \, \left( \bm{\Gamma}_k^{(n)} \right)_{i,j} \, .
    \label{eq:sig_gamma_bar}
\end{eqnarray*}
We further define $\bar{\sigma}^{(n)}$ and $\bar{\gamma}^{(n)}$ as follows:
\begin{equation*}
    \bar{\sigma}^{(n)} \equiv \frac{1}{2m} \sum_{k = 1}^m \, \left[\bar{\sigma}^{(n)}_{2k} + \bar{\sigma}^{(n)}_{2k-1}\right] \;\;\;\; \text{and} \;\;\;\;
    \bar{\gamma}^{(n)} \equiv \frac{1}{m} \sum_{k = 1}^m \, \bar{\gamma}^{(n)}_k \, .
    \label{eq:sigma_gamma_avg}
\end{equation*}
Based on the earlier defined matrices, as for the MOAV, we can also define different quantities according to the matrix of reference and the lags between observations. More specifically, let us first consider the case in which we are interested in lags $h$ such that $0 \leq h < n$. The observations at these lags can belong to the sets of observations within the matrix $\bm{\Sigma}_{2k}^{(n)}$, $\bm{\Sigma}_{2k-1}^{(n)}$ and $\bm{\Gamma}_k^{(n)}$ and, for these sets of observations within matrices $\bm{\Sigma}_{2k}^{(n)}$ and $\bm{\Sigma}_{2k-1}^{(n)}$, we can define the following quantity
\begin{equation*}
    \widetilde{\gamma}(h) \equiv \frac{1}{2m(n-h)} \sum_{k = 1}^{2m} \sum_{s=1}^{n-h}\cov\left(X_{(k-1)n+s}, \,X_{(k-1)n+s+h}\right).
    \label{eq:gamma_avg_1}
\end{equation*}
If, however, the observations at the considered lags are among the set of observations only within the matrix $\bm{\Gamma}_k^{(n)}$, we define the quantity below
\begin{equation*}
    \widetilde{\gamma}^{\ast}(h) \equiv \frac{1}{mh} \sum_{k = 1}^{m} \sum_{s=1}^{h}\cov \left(X_{(2k-1)n+s-h}, \,X_{(2k-1)n+s}\right).
    \label{eq:gamma_avg_2}
\end{equation*}
Finally, when considering lags $h$ such that $n \leq h \leq 2n-1$, the set of observations at these lags can only be considered within the matrix $\bm{\Gamma}_t^{(n)}$ and, for this final case, we define the quantity
\begin{equation*}
    \widetilde{\gamma}(h) \equiv \frac{1}{m(2n-h)} \sum_{k = 1}^{m} \sum_{s=1}^{2n-h}\cov \left(X_{(k-1)n+s}, \,X_{(k-1)n+s+h}\right).
    \label{eq:gamma_avg_3}
\end{equation*}
As for the MOAV case, the above definitions can be seen as generalized definitions of the autocovariance which consists in the average autocovariance for a given lag $h$.
Using the above notations and definitions, we can provide the following result.
\begin{customLemma}{2}
\label{thm:AVNS_NO}
  We define the Non-stationary NOAV as
  \begin{eqnarray*}
      \AV_n &= \frac{1}{2 m n^2} \Bigg\{ 2mn \, \widetilde{\gamma}(0) + 2\Bigg[\sum_{h=1}^{n-1}2m(n-h) \, \widetilde{\gamma}(h)\\
      &- mh \, \widetilde{\gamma}^{\ast}(h) - \sum_{h=n}^{2n-1} m(2n-h) \, \widetilde{\gamma}(h) \Bigg] \Bigg\}.
  \end{eqnarray*}

\end{customLemma}

\begin{IEEEproof}[Proof of Lemma \ref{thm:AVNS_NO}]

The proof of Lemma \ref{thm:AVNS_NO} is direct from the above definitions. Indeed, we have
\begin{eqnarray*}
        \mathbb{E} \left[ \left(\bar{X}_{2k}^{(n)} - \bar{X}_{2k-1}^{(n)} \right)^2\right] &= \var\left(\bar{X}_{2k}^{(n)} \right) + \var\left(\bar{X}_{2k-1}^{(n)} \right) \\
        &- 2 \cov \left(\bar{X}_{2k-1}^{(n)}, \bar{X}_{2k}^{(n)}\right)\\
        &= \bar{\sigma}^{(n)}_{2k} + \bar{\sigma}^{(n)}_{2k-1} - 2 \bar{\gamma}^{(n)}_k.
\end{eqnarray*}
Then, using Eq. (\ref{eq:AVNS_def}), we obtain
\begin{eqnarray*}
        \AV_n &= \frac{1}{2 m} \sum_{k = 1}^{m} \mathbb{E} \left[ \left(\bar{X}_{2k}^{(n)} - \bar{X}_{2k-1}^{(n)} \right)^2\right]\\
        &= \frac{1}{2 m} \sum_{k = 1}^{m} \bar{\sigma}^{(n)}_{2k} + \bar{\sigma}^{(n)}_{2k-1} - 2 \bar{\gamma}^{(n)}_k\\
        &= \frac{1}{2 m} \left[ 2m \, \bar{\sigma}^{(n)} - 2m \,\bar{\gamma}^{(n)} \right]\\
        &= \frac{1}{2 m n^2} \Bigg\{ 2mn \, \widetilde{\gamma}(0) + 2\Bigg[\sum_{h=1}^{n-1}2m(n-h) \, \widetilde{\gamma}(h)\\
        &- mh \, \widetilde{\gamma}^{\ast}(h) - \sum_{h=n}^{2n-1} m(2n-h) \, \widetilde{\gamma}(h) \Bigg] \Bigg\},
\end{eqnarray*}
which concludes the proof.
\end{IEEEproof}

%%%%%%%%%%%%%%%%%%%%%%%%%%%%%%%%%%%%%%%%%%%%%%%%%%%%%%%%%%%%%
%%%%%%%%%%%%%%%%%%%% PROOF THEOREM 1   %%%%%%%%%%%%%%%%%%%%%%
%%%%%%%%%%%%%%%%%%%%%%%%%%%%%%%%%%%%%%%%%%%%%%%%%%%%%%%%%%%%%
\section{}
\label{app.moavns}

\begin{IEEEproof}[Proof of Lemma \ref{lem:AVNS_MO}]
    In order to prove Lemma \ref{lem:AVNS_MO} let $\bar{\sigma}^{(n)}_t$ and $\bar{\gamma}^{(n)}_t$ denote the averages of the matrices $\bm{\Sigma}_{t}^{(n)}$ and $\bm{\Gamma}_t^{(n)}$, respectively, i.e.
\begin{eqnarray*}
    \bar{\sigma}^{(n)}_t \equiv \frac{1}{n^2}\sum_{i = 1}^{n} \sum_{j = 1}^{n} \, \left( \bm{\Sigma}_t^{(n)} \right)_{i,j}, \\
    \bar{\gamma}^{(n)}_t \equiv \frac{1}{n^2} \sum_{i = 1}^{n} \sum_{j = 1}^{n} \, \left( \bm{\Gamma}_t^{(n)} \right)_{i,j} \, .
    \label{eq:sig_gamma_bar}
\end{eqnarray*}
Further, we define $\bar{\sigma}^{(n)}$ and $\bar{\gamma}^{(n)}$ as follows:
\begin{eqnarray*}
    \bar{\sigma}^{(n)} \equiv \frac{1}{2(T-2n+1)} \sum_{t = 2n}^{T} \, \bar{\sigma}^{(n)}_{t} + \bar{\sigma}^{(n)}_{t-n} \\
    \bar{\gamma}^{(n)} \equiv \frac{1}{T-2n+1} \sum_{t = 2n}^{T} \, \bar{\gamma}^{(n)}_t \, .
    \label{eq:sigma_gamma_avg}
\end{eqnarray*}
Based on these definitions we have
\begin{eqnarray*}
        \mathbb{E} \left[ \left(\bar{X}_{k}^{(n)} - \bar{X}_{k-n}^{(n)} \right)^2\right] &= \var\left(\bar{X}_{k}^{(n)} \right) + \var\left(\bar{X}_{k-n}^{(n)} \right)\\
        & - 2 \cov \left(\bar{X}_{k}^{(n)}, \bar{X}_{k-n}^{(n)}\right)\\
        &= \bar{\sigma}^{(n)}_{k} + \bar{\sigma}^{(n)}_{k-n} - 2 \bar{\gamma}^{(n)}_k.
\end{eqnarray*}
Then, using Eq. (\ref{eq:MOAVNS_def}), we obtain
\begin{eqnarray*}
        \AV_n &= \frac{1}{2 m^{\star}} \sum_{k = 2n}^{T} \mathbb{E} \left[ \left(\bar{X}_{k}^{(n)} - \bar{X}_{k-n}^{(n)} \right)^2\right]\\
        &= \frac{1}{2 m^{\star}} \sum_{k = 2n}^{T} \bar{\sigma}^{(n)}_{k} + \bar{\sigma}^{(n)}_{k-n} - 2 \bar{\gamma}^{(n)}_k\\
        &= \frac{1}{2 m^{\ast}} \left[ 2m^{\star} \, \bar{\sigma}^{(n)} - 2m^{\star} \,\bar{\gamma}^{(n)} \right]\\
        &= \frac{1}{2 m^{\star} n^2} \Bigg\{ 2nm^{\star} \, \widetilde{\gamma}(0) + 2\Bigg[\sum_{h=1}^{n-1}2m^{\star}(n-h) \, \widetilde{\gamma}(h)\\
        &- m^{\star}h \, \widetilde{\gamma}^{\star}(h) - \sum_{h=n}^{2n-1} m^{\star}(2n-h) \, \widetilde{\gamma}(h) \Bigg] \Bigg\},
\end{eqnarray*}
which concludes the proof.
\end{IEEEproof}

\end{document}